\documentclass[12pt]{article} 
 
\usepackage{latexsym} 
\usepackage{ntheorem} 
\usepackage[latin1]{inputenc} 

\usepackage{amsfonts,amsmath,amssymb} 

\usepackage{a4wide}

\author{Vincent Bosser (Caen)\footnote{Supported by the contract ANR ``HAMOT", BLAN-0115-01.} \\ 
Andrea Surroca (Basle)\footnote{Supported by an Ambizione fund PZ00P2\_121962 of the Swiss National Science Foundation and the Marie Curie IEF 025499 of the European Community.}}

\title{\Large Upper bound for the height of $S$-integral points on elliptic curves}

\newtheorem{thm}{\textbf{Theorem}}[section]

\newtheorem{lemma}[thm]{\textbf{Lemma}} 
\newtheorem{prop}[thm]{\textbf{Proposition}} 
 
\newtheorem{corollary}[thm]{Corollary}

{\theorembodyfont{\upshape}
}

\newcommand\rat{\mathbf{Q}} 
\newcommand\Q{\mathbf{Q}} 
\newcommand\C{\mathbf{C}} 
\newcommand\Z{\mathbf{Z}} 
\newcommand\enteros{\mathbf{Z}} 
\newcommand\PP{\mathbf{P}} 
\newcommand\R{\mathbf{R}}

\newcommand\card{\mathrm{card}}

\newcommand\ord{\mathrm{ord}}
\newcommand\rk{\mathrm{rk}}
\newcommand\trace{\mathrm{trace}}

\newcommand\Reg{\mathrm{Reg}}

\newcommand\vp{\mathfrak{p}}

 \newcommand\Fu{{K_{0}}}
 \newcommand\Ft{K} 
  
\newcommand\esp{\hspace{0,2cm}}

\input cyracc.def

\begin{document} 
 
\bibliographystyle{alpha}

\maketitle 
 
\begin{quote} 
\textbf{Abstract.} 
{\small   
We establish new upper bounds for the height of the $S$-integral points of an elliptic curve.
This bound is explicitly given in terms of the set $S$ of places of the number field $K$ involved,
but also in terms of the degree of $K$, as well as the rank, the regulator and the height of a basis of the Mordell-Weil group
of the curve.
The proof uses the elliptic analogue of Baker's method, based on lower bounds for linear forms in elliptic logarithms.  
} 

\end{quote}

 

\section{Introduction} 
 
A fundamental problem in Diophantine Geometry is to get effective versions 
of known qualitative results. For example the classical 
finiteness theorem of  Siegel asserts that  the set of integral points of an affine 
algebraic curve of genus greater than one or of genus zero with at least 3 
points at infinity is finite. For that curves of genus greater than 2, Siegel's theorem is superseded by Faltings' theorem which asserts
that the set of rational points is finite. 
These are qualitative statements, but not effective.
To effectively find these points, say, in a fixed number field, it would suffice to find 
an effective upper bound for the height of the points. 
Nowadays, the results of this kind which are known come from Baker's method (based on non trivial lower bounds for linear forms in logarithms). They all concern integral points. The method can be applied for certain classes of curves, in particular, for elliptic curves (\cite{baker-coates}). Generalizing an idea of Gel'fond, S. Lang \cite{lang.dio.an} has shown that 
one can also bound the height of integral points of an elliptic curve $E$ using lower bounds for linear forms in \textit{elliptic} logarithms, in a more natural way than using classical logarithms. 

Let $E$ be an elliptic curve defined over a number field $K_0$, let $K/K_0$ be any finite extension and let
$S$ be a finite set of finite places of $K$. 
In this paper we obtain new upper bounds for the height of the $S$-integral points of $E(K)$, using lower bounds in linear forms in elliptic
logarithms (Theorem \ref{borne-hauteur}).

This method was first applied successfully by D. Masser 
\cite[Appendix IV]{davidmasser} when $\Ft=\Fu=\Q$, the curve $E$ has complex 
multiplication and  $S = \emptyset$. To this end he used his own lower bounds 
for usual (archimedean) elliptic logarithms. D. Bertrand \cite{bertrand} 
then established such lower bounds for $\mathfrak p$-adic elliptic 
logarithms, which allowed him to treat the case $\Ft=\Fu$ and $S$ arbitrary 
(again for curves with complex multiplication). 
Applying the explicit lower bounds for linear forms in elliptic 
logarithms of S. David \cite{david} in the archimedean case, and of 
N. Hirata \cite{noriko} 
in the ultrametric one, we deal here with the general case of an arbitrary 
elliptic curve defined over 
$\Fu$ and of an arbitrary field extension $\Ft/\Fu$. Our results improve the previous results of D. Bertrand. 
Moreover, contrary to the previous works, the bound we obtain for the height 
of the $S$-integral points is not only given in terms of the set $S$, but 
also in terms of the number field $\Ft$. 
More precisely, the ``constant" which occurred in the previous works 
is here explicitly given in terms of the degree $[\Ft:\Q]$, the rank of the Mordell-Weil group $E(\Ft)$, the heights of generators of the
free part of $E(\Ft)$, and the regulator of $E/\Ft$ (but we do not make explicit the
dependence on $E/\Fu$). As mentioned at the end of Section~\ref{upperbound}, 
it is possible to derive from our main result a conditional upper
bound in terms only of the degree $[\Ft:\Q]$, the discriminant of $\Ft$ and the set of places $S$.

\medskip 
 
For convenience to the reader, we have gathered in Section~\ref{notations} the notations which will be used throughout the text. 
We state the main theorem in Section~\ref{section-main} and prove it in Section~\ref{section-bound}.
 
\section{Notations}\label{notations}

Throughout the text, if $x$ is a non negative real number, we set $\log^+x=\max\{1,\log x\}$ (with the convention $\log^+0=1$).
  
If $K$ is a number field, we will denote by $O_K$ its ring of integers, by
$D_K$ the absolute value of its discriminant, and by $M_K$ the set of places
of $K$. The set of the archimedean places will be denoted by $M_K^{\infty}$
and the set of the ultrametric ones will be denoted by $M_K^0$. For each $v$
in $M_K$, we define an absolute value $|\cdot|_v$ on $K$ as follows.
If $v$ is archimedean, then $v$ corresponds to an embedding 
$\sigma : K \hookrightarrow \C$ (we will often identify the place $v$ with
the embedding $\sigma$), and we set $|x|_v=|x|_{\sigma}:=|\sigma(x)|$, where
$|\cdot |$ is the usual absolute value on $\C$. If $v$ is ultrametric, then
$v$ corresponds to a non zero prime ideal $\mathfrak p$ of $O_K$ (we
will identify $v$ and ${\mathfrak p}$),
and we take for
$|\cdot|_v=|\cdot|_{\mathfrak p}$ the absolute value on $K$ normalized
by $|p|_v=p^{-1}$, where $p$ is the prime number such that
${\mathfrak p}\mid p$.
We denote by $K_v$ the completion of $K$ at $v$ and use
again the notation $|\cdot|_v$ for the unique extension of $|\cdot|_v$ to
$K_v$. If $v$ is an ultrametric place associated to the prime ideal
${\mathfrak p}$, we denote by $e_{\mathfrak p}$ the ramification index of
${\mathfrak p}$ over $p$, by $f_{\mathfrak p}$ the residue class degree,
and by $\ord_{\mathfrak p} : K _{\mathfrak{p}}^*\rightarrow\Z$ the valuation
normalized by $\ord_{\mathfrak p}(p)=e_{\mathfrak p}$
(hence $\ord_{\mathfrak p}(x)=-e_{\mathfrak p}\log_p|x|_{\mathfrak p}$ for
all $x$ in $K_{\mathfrak{p}}^*$).
 
If $S$ is a finite subset of $M_K^0$, we denote by 
$$O_{K,S} = \{x \in K; \forall v \notin S\cup M_K^{\infty}, 
|x|_{v} \leq 1  \}$$ 
the ring of $S$-integers of $K$, and we set 
$$\Sigma_S =\sum_{{\mathfrak p}\in S}\log N_{\Ft/\rat}({\mathfrak p}).$$ 
Note that with our notation, the set $S$ contains only non-archimedean places of $K$.
 
Throughout the text, we denote by $h$ the absolute logarithmic Weil height
on the projective space $\PP^n(\overline{\rat})$, and we denote by
$h_K:=[K:\rat]h$ the relative height on $\PP^n(K)$. Thus,
if $(\alpha_{0}: \ldots: \alpha_{n}) \in \mathbf{P}^{n}(K)$, we have:
\begin{equation}\label{h.weil} 
h(\alpha_{0}: \ldots : \alpha_{n}) = \frac{1}{[K:\rat]} 
\sum_{v \in M_{K}} [K_v:\Q_v] \log \max \{|\alpha_{0}|_{v}, \ldots, |\alpha_{n}|_{v}\}. 
\end{equation} 

\medskip 
 
Let $E\subset\PP^2$ be an elliptic curve defined over a number field $K$.
The Mordell-Weil group  $E(K)$ of $K$-rational points of $E$ is a finitely generated group:
$$E(K) \simeq E(K)_{tors} \oplus \enteros^{\rk(E(K))}.$$ 
We will often simply write  $r = \rk(E(K))$ for its rank, and we will
denote by $(Q_1,\ldots,Q_r)$ a basis of its free part. We will also
denote by $O$ the zero element of $E(K)$.

We  further denote by $\hat{h}:E(\overline{K})\rightarrow\R$ the N\'eron-Tate height on $E$. 
The ``N\'eron-Tate pairing'' $<\, ,\, >$ is defined by $<P,Q> = \frac{1}{2}(\hat{h}(P+Q) - \hat{h}(P) - \hat{h}(Q))$.
The regulator $\Reg(E/K)$ of $E/K$ is the determinant of the matrix $\mathcal{H}=(<Q_{i}, Q_{j}>)_{1\leq i, j\leq r}$ of the N\'eron-Tate pairing with respect to the chosen basis $(Q_1,\ldots,Q_r)$, that is 
$$\Reg(E/K) =  \det(\mathcal{H}).$$

If the elliptic curve is defined by a Weiertrass equation
$y^{2} = x^{3} + Ax + B$ with $A, B$ in $O_K$, then we have the origin $O=(0:1:0)$.
If $Q\ne O$ is a point of $E$, we then denote its affine coordinates
(in the above Weierstrass model) as usual by $(x(Q),y(Q))$.
For $Q$ in $E(\overline{\Ft})$ we define $h_x(Q):=h(1:x(Q))$ if $Q\not=O$ and $h_x(O):=0$. 
Finally, we denote by $E(O_{K,S})$ the set of $S$-integral points of $E(K)$
with respect to the $x$-coordinate, that is
$$E(O_{K,S}) = \{ Q \in E(K)\setminus\{O\}; 
x(Q) \in O_{K,S}  \} \cup \{O\}.$$

In the whole text, we will fix a number field $\Fu$ and an elliptic curve $E$ defined over $\Fu$.
Since we do not explicit any dependence on $E/\Fu$, we will call ``constant" any quantity depending on $E/\Fu$.
This convention about constants will apply in particular to the implicit constant involved in the symbol $\ll$,
where $X \ll Y$ means here that $X \leq c(E/\Fu) Y$, where $c(E/\Fu) \geq 1$ is a number depending at most on $E/\Fu$.

\section{Statement of the result}\label{section-main}

Let $\Fu$ be a fixed number field, and let $E\subset\PP^2$ be
an elliptic curve defined by a Weierstrass equation
\begin{equation}\label{weierstrass} 
y^{2} = x^{3} + Ax + B 
\end{equation} 
with $A,B\in O_{\Fu}$.
Let  $\Ft$  be a finite extension of $\Fu$ and $S\subset M_{\Ft}^0$ a finite set of places of $K$.

According to the notations of Section~\ref{notations}, we
put $r=\rk(E(K))$, we denote by $(Q_1,\ldots,Q_r)$ any basis of the free part of the Mordell-Weil group $E(K)$, and we write $\Reg(E/K)$ for the regulator of $E/K$.
We further set 
$$d:=[\Ft:\Q],$$
 and we define the real number $V$ by
$$\log V:=\max\{\hat{h}(Q_i); 1\le i\le r\}.$$

\smallskip 
 
The main result of this article is the following: 
 
\begin{thm}\label{borne-hauteur}  In the above set up, let $Q$ be a point in $E(O_{\Ft,S})$. Then there exist positive effectively computable real numbers
$\gamma_0, \gamma_1$ and $\gamma_2$ depending only on $A$ and $B$ (that is, on the curve $E/\Fu$), such that,
if $r = 0$, then $h_x(Q) \le \gamma_0$, and, if $r > 0$, then
\begin{equation}\label{bound} 
h_x(Q) \le C_{E,\Ft}e^{(8r^2+\gamma_1dr)\Sigma_S},
\end{equation} 
where 
\begin{align}\label{bound2}
C_{E,\Ft} & = \gamma_2^{r^2}r^{2r^2} d^{9r+15}(\log^+d)^{r+6} 
(\log^+\log V)^{r+7}  (\log^+\log^+\log V)^2\, \prod_{i=1}^{r}\max\{1,\hat h(Q_i)\} \cr 
& \times \log^+\!(\Reg(E/K)^{-1}) 
(\log^+\!\log(\Reg(E/K)^{-1}))^2 
(\log^+\!\log^+\!\log(\Reg(E/K)^{-1})).
\end{align}  
\end{thm}

The bounds obtained by classical Baker's method often depend on $d$, $D_K$ and $\Sigma_S$.
The bound of Theorem \ref{borne-hauteur} depends on $d$ and $\Sigma_S$, as well as on the rank $\rk (E(K))$, the heights $\hat{h}(Q_i)$
and the regulator of $E/K$.
Because of the different nature of the parameters involving $K$, it makes sense to compare our result with the results obtained by classical Baker's method
only when the number field is fixed. 

Denote by $s$ the cardinal of $S$, by $P(S):=\{p\textrm{ prime}\mid \exists v\in S, v|p\}$ the residue characteristics of $S$, and by $P$
its maximum.  For $K= \rat$, L. Hajdu and T. Herendi \cite{hajdu-herendi} obtained the following result:   
$$\max \{h(x), h(y)\} \leq (\kappa_{1} \,s + \kappa_{2})\, 10^{38\, s + 86}\, (s + 1)^{20\,s + 35} \,P^{24}\, (\log^{+}(P))^{4\,s +2},$$
where $\kappa_{1}$ and $\kappa_{2}$ depend at most on $A$ and $B$. For any  $K$,  the Corollary 6.9 of \cite{bornes} gives: 
$$h_{x}(Q)\leq k_0 c_{d}^{s + k_{1}}\,s^{20\,s + k_{2}}\,P^{4\,d}\,(\log P)^{8\,s + k_{3}}\, e^{\gamma({E,K,S})}\,\gamma({E,K,S})^{8d-2},$$
where the numbers $k_{i}$ depend only on  $E/\Fu$, $c_d$ depends only on the degree $d$ and $\gamma({E,K,S}) = 4\,(\log D_{K} + \Sigma_{S} + k_{4} + 4\, \log 4 \, d\,\frac{\Sigma_{S} + k_{5}}{\log(\Sigma_{S} + k_{5})})$.
In order to compare these different bounds we may use the following inequalities (for the last one, one may use the Prime Number Theorem, see, for example \cite[Lemma 2.1]{bornes}):
$$(d \cdot \card P(S))^{-1} \Sigma_S \leq \log P \leq \Sigma_S, \textrm{and}$$
$$\card (S) \leq d \cdot \card P(S) \leq 4 d \frac{\Sigma_S}{\log \Sigma_S}.
$$
One can see that, for a fixed $K$, the bounds of \cite[Corollary 6.9 ]{bornes} and Theorem \ref{borne-hauteur} are of the same order
and that for $K=\rat$ the bound of \cite{hajdu-herendi} is stronger.

\section{Proof of Theorem~\ref{borne-hauteur}}\label{section-bound} 

To prove Theorem~\ref{borne-hauteur}, we will 
determine an upper bound for $|x(Q)|_v$ for each $v\in M_{\Ft}$ and then 
sum over all the places $v$ using the formula (\ref{h.weil}). 
The upper bound for $|x(Q)|_v$ will be obtained using the explicit 
lower bounds for linear forms in elliptic logarithms of 
\cite{noriko} and \cite{david}. 
In the next section, we first treat the case of an archimedean place $v$. 
Then, in Section~\ref{nonarch}, we handle the case where $v$ is ultrametric. 
We can then prove Theorem~\ref{borne-hauteur} in Section~\ref{upperbound}.

In the next sections, we denote by $\kappa_1, \kappa_{2}, \ldots, c_1, c_{2}, \ldots$ positive real numbers
(which we will call ``constants") depending at most on $A$ and $B$, \textit{i.e.} on the curve $E/\Fu$. We use the greek letters for the constants appearing in the statements and the latin letters for the proofs. In each proof we start counting by $c_1$. 
  
\subsection{The archimedean case}\label{arch}

In this section we fix an archimedean place $v$ of $\Ft$, and we
assume that the rank $r$ of the group $E(\Ft)$ is non zero.
We denote by $\sigma:\Ft\hookrightarrow \C$ the embedding corresponding
to $v$, and by $E_{\sigma}$ the elliptic curve defined by
$$y^2=x^3+\sigma(A)x+\sigma(B).$$
The homomorphism $\sigma$ obviously induces a group isomorphism
$\sigma:E(\Ft)\simeq E_{\sigma}(\sigma(\Ft))$.
 
Put $g_{2,\sigma}=-4\sigma(A)$ and $g_{3,\sigma}=-4\sigma(B)$. Then
$E_{\sigma}$ is isomorphic to the elliptic curve defined by
\begin{equation}
Y^2=4X^3-g_{2,\sigma}X-g_{3,\sigma}
\end{equation}
under the substitution $X=x$, $Y=2y$.
Let $\Lambda_{\sigma}$ be the lattice of $\C$ with invariants
$g_{2,\sigma}$, $g_{3,\sigma}$. We will consider
the exponential map of $E_{\sigma}$, which is given by
\begin{equation}
\begin{array}{cccl}
\exp_{\sigma}: & \C & \rightarrow & E_{\sigma}(\C)\cr 
& z & \mapsto & \left\{\begin{array}{l} 
(\wp_{\sigma}(z):\wp_{\sigma}'(z)/2:1)\ \text{if}\ z\notin\Lambda_{\sigma}\cr 
                (0:1:0)\ \text{if}\ z\in\Lambda_{\sigma}, \end{array}\right. 
\end{array}
\end{equation}
where $\wp_{\sigma}$ is the Weierstrass function associated to the lattice
$\Lambda_{\sigma}$.
It induces a group isomorphism $\C/\Lambda_{\sigma}\simeq E_{\sigma}(\C)$.
Let $(\omega_{1,\sigma},\omega_{2,\sigma})$ be a basis of the lattice
$\Lambda_{\sigma}$,
and denote by $\Pi_{\sigma}$ the associated fundamental parallelogram
centered at zero.
Then $\exp_{{\sigma}}|_{\Pi_{\sigma}}:\Pi_{\sigma}\rightarrow E_{\sigma}(\C)$
is bijective and we will denote
by $\psi_{\sigma}: E_{\sigma}(\C)\rightarrow\Pi_{\sigma}$ its inverse map
(the elliptic logarithm).

\begin{prop}\label{archi}
Let $Q$ be a non-torsion point of $E(\Ft)$. Write
$Q=m_1Q_1+\cdots +m_rQ_r+Q_{r+1}$, where $Q_{r+1}\in E(\Ft)$ is a torsion
point and $m_i\in\Z$, $1\le i\le r$, and define
$M:=\max\{|m_1|,\ldots,|m_r|\}$. Recall that $\log V:=\max\{\hat{h}(Q_i); 1\le i\le r\}$. Then we have
\begin{align*}
\log |x(Q)|_{\sigma}\le
\kappa_{1}^{r^2}r^{2r^2} d^{2r+8} (\log^+d)^{r+5} &(\log^+ M)(\log^+\log M)\cr
&\times (\log^+\log V)^{r+5}
\prod_{i=1}^{r}\max\{1,\hat h(Q_i)\}.
\end{align*}
\end{prop}
 
To prove this we will use the following result, which is a consequence of
Theorem~2.1 of \cite{david} (in which we have chosen $E=e$).

\begin{thm}[S. David]\label{sinnou}
Let $m_1,\ldots,m_{r+1},n_1,n_2$ be rational integers,
and let $\gamma_1,\ldots,\gamma_{r+1}$ be elements of
$E_{\sigma}(\sigma(\Ft))$.
Define $u_i=\psi_{\sigma}(\gamma_i)$, $1\le i\le r+1$, and put 
$$\mathcal{L}=m_1u_1+\cdots+m_{r+1}u_{r+1}+n_1\omega_{1,\sigma}+n_2 
\omega_{2,\sigma}.$$ 
Set further $B=\max\{|m_1|,\ldots,|m_{r+1}|,|n_1|,|n_2|\}$ and $\log W= 
\max\{\hat h(\gamma_i),1\le i\le r+1\}$. 
If $\mathcal{L}\not=0$, then 
\begin{align*}  
\log |\mathcal{L}| \ge 
-\kappa_{2}^{r^2}r^{2r^2} d^{2r+8}(\log^+d)^{r+5}&(\log^+B)(\log^+\log B)\cr 
&\times (\log^+\log W)^{r+5} \prod_{i=1}^{r+1}\max\{1,\hat h(\gamma_i)\}. 
\end{align*} 
\end{thm} 
 
\noindent {\em Proof of Proposition~\ref{archi}.} 
Let $Q$, $Q_{r+1}$, $m_1,\ldots,m_r$ and $M$ be as in the Proposition.
We set $\gamma_i=\sigma(Q_i)\in E_{\sigma}(\C)$ and
$u_i=\psi_{\sigma}(\gamma_i)$, $1\le i\le r+1$.
We have $\sigma(Q)=m_1\gamma_1+\cdots+ m_r\gamma_r+
\gamma_{r+1}$ and
thus, by definition of $\psi_{\sigma}$ there exist $n_1,n_2\in\Z$ such that
\begin{equation}\label{logell}
\psi_{\sigma}(\sigma(Q))= m_1u_1+\cdots +m_ru_r+u_{r+1}+n_1\omega_{1,\sigma}+
n_2\omega_{2,\sigma}.
\end{equation}
Since the function $\wp_{\sigma}$ has a pole of order $2$ at zero,
there exists a constant $c_{1,\sigma}=c_{1}(\sigma(A),\sigma(B))>0$ such
that 
$|z^2\wp_{\sigma}(z)|\le c_{1,\sigma}$ for all $z$ in $\Pi_{\sigma}$.
Applying this to $z=\psi_{\sigma}(\sigma(Q))$ and
putting $c_{1}:=\max_{\sigma}\{c_{1,\sigma}\}$ (which depends only on the restriction of $\sigma$ to $\Fu$ hence on $E/\Fu$ only),
we find that
\begin{equation*}
|x(Q)|_{\sigma}\le c_{1} \cdot |\psi_{\sigma}(\sigma(Q))|^{-2}
\end{equation*}
i.e.
\begin{equation}\label{www}
\log |x(Q)|_{\sigma}\le \log c_{1}-2\log |\psi_{\sigma}(\sigma(Q))|.
\end{equation}
In order to use Theorem~\ref{sinnou}, observe that since all the norms on $\R^2$ are equivalent, there exists a
constant $c_{2,\sigma}=c_{2}(\sigma(A),\sigma(B))>0$ such that
$|x\omega_{1,\sigma}+y\omega_{2,\sigma}|\ge c_{2,\sigma}
\max\{|x|,|y|\}$ for all real numbers $x,y\in\R$. Therefore
we have, using (\ref{logell}) and since obviously
$|u|\le (|\omega_{1,\sigma}|+|\omega_{2,\sigma}|)/2$ for every
$u$ belonging to the fundamental parallelogram $\Pi_{\sigma}$,
$$c_{2,\sigma}\max\{|n_1|,|n_2|\}\le
|\psi_{\sigma}(\sigma(Q))|+ M(|u_1|+\cdots +|u_{r+1}|)\le
c_{3,\sigma}(1+(r+1)M)\leq c_{4,\sigma}r M.
$$
Hence
$$ \max\{|n_1|,|n_2|\} \leq c_{2} r M$$
with $c_{2}=\max_{\sigma}\{c_{4,\sigma}/c_{2,\sigma}\}$
(which depends only on $E/\Fu$ for the same reason as above), and so
\begin{equation}\label{coeffs}
B:= \max\{1, |m_1|,\ldots,|m_r|,|n_1|,|n_2|\}\le c_{3} r M.
\end{equation}
Applying now Theorem~\ref{sinnou} to the linear form
(\ref{logell}) (which is not zero since $Q$ is not torsion) and taking into account (\ref{coeffs}), we deduce
\begin{align*}  
\log |\psi_{\sigma}(\sigma(Q))| \ge 
-c_{4}^{r^2}r^{2r^2} d^{2r+8}(\log^+d)^{r+5}&(\log^+M)(\log^+\log M)\cr 
&\times (\log^+\log V)^{r+5} \prod_{i=1}^{r}\max\{1,\hat h(Q_i)\}. 
\end{align*} 
This estimate, together with (\ref{www}), yields the proposition.
\hfill$\Box$

\subsection{The ultrametric case}\label{nonarch}
 
We fix here an
ultrametric place $v$ of $\Ft$ associated to a prime ideal $\vp$
lying above the prime number $p$, and we assume again that the rank $r$
of the group $E(\Ft)$ is non zero. We will prove :
 
\begin{prop}\label{ultra} 
Let $Q$ be a non-torsion point of $E(\Ft)$.
Write $Q=m_1Q_1+\cdots +m_rQ_r+Q_{r+1}$, where $Q_{r+1}\in E(\Ft)$ is a torsion
point and $m_i\in\Z$, $1\le i\le r$, and define
$M:=\max\{|m_1|,\ldots,|m_r|\}$. Recall that $\log V:=\max\{\hat{h}(Q_i); 1\le i\le r\}$. Then we have
\begin{align*}  
\log  |x(Q)|_{\vp}\le 
\kappa_3^{r^2} r^{2r^2} p^{8r^2+\kappa_4dr} d^{9r+14}&(\log^+d)^{r+3} (\log^+ M)\cr 
& \times (\log^+\log  V)^{r+3}\prod_{i=1}^{r}\max\{1,\hat h(Q_i)\}.
\end{align*} 
\end{prop} 

To prove this Proposition we will use the $v$-adic 
exponential map of $E$, whose definition and properties we now recall 
for reader's convenience. 
By \cite{weil} (see also \cite{lutz}), there exists a unique function $\psi(z)$ 
analytic in a neighbourhood of $0$ in ${\Ft}_{\mathfrak p}$ which satisfies 
$$ \psi'(z)=(1+A z^4 + Bz^6)^{-\frac{1}{2}}; \esp \psi(0)=0$$
(where we define of course $(1+t)^{-1/2}=1-t/2+\cdots$ for
$t$ in ${\Ft}_{\mathfrak p}$ with $|t|_{\mathfrak p}$ small).
It is not difficult to see that this function $\psi$ is analytic in the open disk
$$ \mathcal{C}_{\vp} = \{z \in {\Ft}_{\vp}; |z|_{\vp} < p^{-\lambda_p} \},$$
where $\lambda_p = (p-1) ^{-1}$ if $p\not= 2$ and $\lambda_p =1/2$ if $p=2$.
Moreover, one can show that for $z\in\mathcal{C}_{\vp}$, $z\not=0$, 
we have $$ \psi(z)= z+\sum_{n\ge 2}\psi_nz^n$$ 
with $|\psi_nz^n|_{\vp}<|z|_{\vp}$ for all $n\ge 2$. It follows from 
results of non-archimedean analysis (see {\em e.g.} \cite[Satz~2]{Gun}) 
that $\psi$ induces a bijection 
$\psi:\mathcal{C}_{\vp}\rightarrow\mathcal{C}_{\vp}$, whose inverse map 
is also analytic. Let $\varphi=\psi^{-1}:\mathcal{C}_{\vp}\rightarrow 
\mathcal{C}_{\vp}$ be this inverse map. Then $\varphi$ is the unique 
solution on $\mathcal{C}_{\vp}$ of the differential equation 
$$y' = (1+A y^4 + By^6)^{\frac{1}{2}}; \esp y(0) =0.$$ 
Moreover, $|\varphi(z)|_{\vp}=|z|_{\vp}$ and $\varphi(-z)=-\varphi(z)$ 
for all $z\in\mathcal{C}_{\vp}$ 
(note that the similar results proved in \cite[p.~246]{lutz} give 
a smaller disk than our $\mathcal{C}_{\vp}$ when $p=2$). 
 
Set now $\wp=1/\varphi^2$. One has on $\mathcal{C}_{\vp}\setminus\{0\}$ 
$$ \frac{1}{4}\wp'^2=\wp^3+A\wp+B.$$ 
The $v$-adic exponential map of $E$ is then defined by 
$$\begin{array}{rccl} 
\exp_{\vp} : & \mathcal{C}_{\vp} & \rightarrow &E({\Ft}_{\vp})\\ 
       & z             & \mapsto & \left\{\begin{array}{l} 
                                (\wp(z):\wp'(z)/2: 1)\mbox{\ if\ } z\not=0\cr 
                                     O \mbox{\ if\ } z=0\end{array} \right. 
\end{array}$$ 
or equivalently by 
\begin{equation}\label{exponentielle} 
\exp_{\vp}(z)=(\varphi(z):-\varphi'(z): \varphi^3(z)) 
\end{equation} 
for all $z\in\mathcal{C}_{\vp}$. This is an injective group homomorphism 
which is not surjective. Let 
$\mathcal{U}_{\vp} = \exp_{\vp}(\mathcal{C}_{\vp})$ be the image of the 
exponential map. It is known that the group 
$E({\Ft}_{\vp})/\mathcal{U}_{\vp}$ is finite. We will need 
an explicit upper bound for the exponent of this group. 
 
\begin{lemma}\label{exposant} 
The exponent $\nu_{\vp}$ of the group 
$E({\Ft}_{\mathfrak{p}})/{\cal U}_{\mathfrak{p}}$ 
satisfies 
$$\nu_{\vp} \le  p^{\kappa_5d}.$$ 
\end{lemma}

\noindent {\em Proof.} In what follows we will denote by 
${\mathfrak{O}}_{\mathfrak{p}}=\{z\in {\Ft}_{\mathfrak{p}}; 
|z|_{\mathfrak{p}}\le 1\}$ 
the valuation ring of ${\Ft}_{\vp}$, by 
${\mathfrak{M}}_{\mathfrak{p}}=\{z\in {\Ft}_{\mathfrak{p}};  
|z|_{\mathfrak{p}}<1\}$ the maximal ideal of ${\mathfrak{O}}_{\mathfrak{p}}$, 
by $\pi$ a uniformizer ({\em i.e.} ${\mathfrak{M}}_{\mathfrak{p}}=\pi 
{\mathfrak{O}}_{\mathfrak{p}}$), and 
by $k(\vp)={\mathfrak{O}}_{\mathfrak{p}}/{\mathfrak{M}}_{\mathfrak{p}}$ 
the residue field of $\mathfrak{O}_{\vp}$. 
 
Let $E_{\mathfrak{p}}\subset\PP_2$ be a minimal Weierstrass model of 
$E$ at ${\mathfrak{p}}$. We know that there is an admissible change 
of coordinates 
\begin{equation}\label{change} 
f:{\Ft}_{\mathfrak{p}}^2\rightarrow {\Ft}_{\mathfrak{p}}^2,\quad 
(x,y)\mapsto (u^2x+r,u^3y+u^2sx+t) 
\end{equation} 
with $u,r,s,t\in {\mathfrak{O}}_{\mathfrak{p}}$, 
such that $f$ induces a group isomorphism 
$f:E_{\mathfrak{p}}({\Ft}_{\mathfrak{p}})\simeq E({\Ft}_{\mathfrak{p}})$. 
Thus, it suffices to estimate the exponent of the group 
$E_{\mathfrak{p}}({\Ft}_{\mathfrak{p}})/{\mathfrak{U}}_{\mathfrak{p}}$, 
where 
${\mathfrak{U}}_{\mathfrak{p}}:=f^{-1}({\cal U}_{\mathfrak{p}})$. 
 
We now claim that the group ${\cal U}_{\mathfrak{p}}$ is explicitly given by 
\begin{equation}\label{U_p} 
{\cal U}_{\mathfrak{p}} = \{(x:y:1)\in E({\Ft}_{\mathfrak{p}});\ 
|x|_{\mathfrak{p}}>p^{2\lambda_p}\}\cup\{O\}. 
\end{equation} 
Indeed, ${\cal U}_{\mathfrak{p}}$ is clearly contained in the right-hand side 
of (\ref{U_p}) since $|\varphi(z)|_{\mathfrak{p}}=|z|_{\mathfrak{p}}$ 
for all $z\in\mathcal{C}_{\vp}$. Conversely, if 
$(x:y:1)\in E({\Ft}_{\mathfrak{p}})$ satisfies 
$|x|_{\mathfrak{p}}>p^{2\lambda_p}$, then we can write 
$$ y^2/x^2= x(1+\frac{A}{x^2}+\frac{B}{x^3})= x(1+t)$$ 
with 
$$ |t|_{\mathfrak{p}} = |\frac{A}{x^2}+\frac{B}{x^3}|_{\mathfrak{p}} 
<p^{-4\lambda_p}.$$ 
It follows that $(1+t)$ is a square in ${\Ft}_{\mathfrak{p}}$ since then 
the series $(1+t)^{1/2}=1+t/2-t^2/8+\cdots$ converges 
in ${\Ft}_{\mathfrak{p}}$ (for $p=2$ it converges as soon as 
$|t|_{\mathfrak{p}}<|2|^2_{\mathfrak{p}}$), and thus 
$x$ is also a square in ${\Ft}_{\mathfrak{p}}$, say $x=\alpha^2$. 
We then have $\alpha^{-1}$ in $\mathcal{C}_{\vp}$, and it follows that 
there exists $z$ in $\mathcal{C}_{\vp}$ such that $\alpha^{-1}=\varphi(z)$, 
{\em i.e.} $x=\varphi(z)^{-2}$. We have moreover 
$y^2=x^3+Ax+B=\varphi'^2(z)/\varphi^6(z)$. Hence, taking $-z$ 
instead of $z$ if necessary, we may choose $z$ so that 
$y=-\varphi'(z)/\varphi^3(z)$. 
Therefore, we have found $z$ in $\mathcal{C}_{\vp}$ such that 
$\exp_{\mathfrak{p}}(z)=(x:y:1)$. This proves (\ref{U_p}). 
 
Using the formulas (\ref{change}) and the ultrametric inequality, we 
deduce from this 
\begin{align*}
{\mathfrak{U}}_{\mathfrak{p}} & 
= \{(x:y:1)\in E_{\mathfrak{p}}({\Ft}_{\mathfrak{p}});\ |x|_{\mathfrak{p}}> 
|u|_{\mathfrak{p}}^{-2}p^{2\lambda_p}\}\cup\{0\}\\ 
& = \{(x:y:1)\in E_{\mathfrak{p}}({\Ft}_{\mathfrak{p}});\ 
\ord_{\mathfrak{p}}(x)/2< -(e_{\mathfrak{p}}\lambda_p+\ord_{\mathfrak{p}}(u))\} 
\cup\{0\}. 
\end{align*} 
In other words, if we denote the canonical ${\mathfrak{p}}$-adic filtration 
of $E_{\mathfrak{p}}$ as in \cite{Hus}, Chapter~14, by 
$$ E_{\mathfrak{p}}({\Ft}_{\mathfrak{p}})\supset 
E_{\mathfrak{p}}^{(0)}({\Ft}_{\mathfrak{p}})\supset\cdots\supset 
E_{\mathfrak{p}}^{(n)}({\Ft}_{\mathfrak{p}})\supset\cdots, $$ 
we see that ${\mathfrak{U}}_{\mathfrak{p}}= 
E_{\mathfrak{p}}^{(n)}({\Ft}_{\mathfrak{p}})$ 
with $n=[e_{\mathfrak{p}}\lambda_p]+\ord_{\mathfrak{p}}(u)+1$. 
 
Estimating the exponent of the group 
$E_{\mathfrak{p}}({\Ft}_{\mathfrak{p}})/{\mathfrak{U}}_{\mathfrak{p}} 
=E_{\mathfrak{p}}({\Ft}_{\mathfrak{p}})/E_{\mathfrak{p}}^{(n)}(\Ft_{\mathfrak{p}})$ 
now easily follows from well-known properties of the 
${\mathfrak{p}}$-adic filtration. Indeed, let $\Delta_{\vp}\in {\Ft}_{\vp}$ be the minimal 
discriminant of the elliptic curve $E$ at ${\vp}$.
By the addendum to Theorem~3 of \cite{tate.AEC} we first have 
\begin{equation}\label{addendum} 
[E_{\mathfrak{p}}(\Ft_{\mathfrak{p}}):E_{\mathfrak{p}}^{(0)}(\Ft_{\mathfrak{p}})] 
\le \max\{4,\ord_{\mathfrak{p}}(\Delta_{\vp})\}, 
\end{equation} 
and by \cite[Proposition~VII.2.1]{silverman.AEC} we have 
\begin{equation}\label{indice.fini} 
[E_{\mathfrak{p}}^{(0)}(\Ft_{\mathfrak{p}}): 
E_{\mathfrak{p}}^{(1)}(\Ft_{\mathfrak{p}})]\le 2\card(k(\vp))+1= 
2p^{f_{\mathfrak{p}}}+1\le\frac{5}{2}p^{f_{\mathfrak{p}}}. 
\end{equation} 
On the other hand, if we define $\widehat{\mathfrak{M}^m_{\mathfrak{p}}}$ 
for every $m\ge 1$ as 
the set $\mathfrak{M}^m_{\mathfrak{p}}$ endowed with the group structure 
given by the formal group law associated to $E_{\mathfrak{p}}$, 
we know that the map 
$t:E_{\mathfrak{p}}^{(m)}(\Ft_{\mathfrak{p}})\rightarrow 
\widehat{\mathfrak{M}^m_{\mathfrak{p}}}$ 
defined by $t(O)= 0$ and $t(Q)=-x(Q)/y(Q)$ if $Q\not=O$ 
is a group isomorphism (see {\em e.g.} the proof of Theorem~14.1.2 of 
\cite{Hus}). It follows, by \cite[Proposition~IV.3.2(a)]{silverman.AEC}, 
that we have for every $m\ge 1$ group isomorphisms 
$$ E_{\mathfrak{p}}^{(m)}(\Ft_{\mathfrak{p}})/ 
E_{\mathfrak{p}}^{(m+1)}(\Ft_{\mathfrak{p}})\simeq 
\widehat{\mathfrak{M}^m_{\mathfrak{p}}}/ 
\widehat{\mathfrak{M}^{m+1}_{\mathfrak{p}}}\simeq 
\mathfrak{M}^m_{\mathfrak{p}}/\mathfrak{M}^{m+1}_{\mathfrak{p}} 
\simeq k(\vp).$$ 
Since the characteristic of the field $k(\vp)$ is equal to $p$, 
we thus get that the exponent of the group 
$ E_{\mathfrak{p}}^{(m)}(\Ft_{\mathfrak{p}})/ 
E_{\mathfrak{p}}^{(m+1)}(\Ft_{\mathfrak{p}})$ is equal to $p$ for all $m\ge 1$. 
Hence we deduce, using (\ref{addendum}) and (\ref{indice.fini}), 
that the exponent of the group 
$E_{\mathfrak{p}}(\Ft_{\mathfrak{p}})/E_{\mathfrak{p}}^{(n)}(\Ft_{\mathfrak{p}})$ 
is at most 
$$ 
\frac{5}{2}\max\{4,\ord_{\mathfrak{p}}(\Delta_{\vp})\}p^{f_{\mathfrak{p}}+n-1}. 
$$ 
Let $\Delta=-16(4A^3+27B^2)$ be the discriminant of the equation 
(\ref{weierstrass}). 
Write now $\Delta=u^{12}\Delta_{\vp}$.
We have
$$n-1=[e_{\vp}\lambda_p]+(\ord_{\vp}(\Delta) 
-\ord_{\vp}(\Delta_{\vp}))/12\le [e_{\vp}\lambda_p]+\ord_{\vp}(\Delta)/12,$$ 
hence
$$\nu_{\vp}\le\frac{5}{2}\max \{4,\ord_{\mathfrak{p}}(\Delta_{\vp})\} 
p^{f_{\mathfrak{p}}+[e_{\mathfrak{p}}\lambda_p] 
+\ord_{\mathfrak{p}}(\Delta)/12}.$$
Noticing now that $\ord_{\vp}(\Delta_{\vp})\le
\ord_{\vp}(\Delta)\le c_1 e_{\vp}$ (with
$ c_1= \displaystyle{\max_{v\in M^0_{\Q(A,B)}}}\{
\ord_v(\Delta)\}$) and since $\lambda_p\le 1/2$, we find
$$ 
\nu_{\vp}\le c_{2}e_{\vp}p^{f_{\vp}+e_{\vp}/2+c_1e_{\vp}/12}\le c_{2}d 
p^{c_{3}d}.$$
\hfill$\Box$ 
 
\medskip 

We will also need the following lemma (where we set $x(O)=\infty$):

\begin{lemma}\label{filtr}
Let $Q$ be a point of $E(\Ft_{\vp})$ such that $|x(Q)|_{\vp}>1$.
Then, for any positive integer $m$, we have $|x(mQ)|_{\vp}\geq |x(Q)|_{\vp}$.
\end{lemma}

\noindent {\em Proof.} Let
$$ E({\Ft}_{\mathfrak{p}})\supset
E^{(0)}({\Ft}_{\mathfrak{p}})\supset\cdots\supset 
E^{(n)}({\Ft}_{\mathfrak{p}})\supset\cdots $$ 
denote the canonical $\vp$-adic filtration of $E$
(see for instance \cite{Hus}, Section~14.1). We recall
that for $n\ge 1$ we have
\begin{equation}\label{filtration}
E^{(n)}({\Ft}_{\mathfrak{p}}) = \{Q\in E(K_{\vp});
\ord_{\vp}(x(Q))\leq -2n\}.
\end{equation}
Let $Q\not =O$ be a point of $E(\Ft_{\vp})$ as in the lemma. We know that $\ord_{\vp}(x(Q))$ is even and thus $Q$ belongs to
$E^{(n)}(K_{\vp}$) with $n:=-\ord_{\vp}(x(Q))/2\geq 1$. Since
$mQ$ also belongs to $E^{(n)}(K_{\vp})$ ($E^{(n)}(K_{\vp})$ is a
group), it follows at once from (\ref{filtration}) that
$\ord_{\vp}(x(mQ))\leq -2n=\ord_{\vp}(x(Q))$.
\hfill $\Box$

The following Theorem was kindly communicated to us by N. Hirata (see \cite{noriko}) :

\begin{thm}[N. Hirata]\label{hirata}
Let $ \beta_1,\ldots,\beta_{n}$ be elements of $\Ft$,
and let  $\gamma_1,\ldots,\gamma_{n}$ be $n$ elements of
$E(\Ft)\cap {\cal U}_{\mathfrak{p}}$. Define
$u_i=\exp^{-1}_{\vp}(\gamma_i)$, $1\le i\le n$, and let 
$$\mathcal{L}=\beta_1u_1+\cdots+\beta_{n}u_{n}.$$
Define the following parameters :
$$\log B =\max\{1,h(\beta_1),\ldots,h(\beta_{n})\}$$
$$ h_E= \max\{1, h(1:A:B)\}$$
$$ \mathcal{E} = \frac{p^{-\lambda_p}}{\displaystyle\max_{1\le i\le n}\{|u_i|_{\vp}\}}$$
$$ \delta = \max\{1,\frac{d}{\log \mathcal{E}} \}$$
$$ g = \max_{1\le i\le n}\{1,\log\delta, h_E, \log\hat h(\gamma_i)\}$$
If $\mathcal{L}\not=0$, then
\begin{align}\label{min_hirata}
\log |\mathcal{L}|_{ \vp}\ge
- \kappa_6^{n^2} (n+1)^{2n(n+8)} p^{8n(n+1)}&\delta^{2n+2}(\log
\mathcal{E})^{-2n-1}(\log B+g+\log(\delta \mathcal{E}))\cr
&\times (g+\log(\delta \mathcal{E}))^{n+1} \prod_{i=1}^{n}(h_E+\max\{1,\hat h(\gamma_i)\}),
\end{align}
where $\kappa_6>0$ is an absolute constant.
\end{thm}

\begin{corollary}\label{min-ultra}
Let $m_1,\ldots,m_{r+1}$ be rational integers, 
and let  $\gamma_1,\ldots,\gamma_{r+1}$ be $r+1$ elements of 
$E(\Ft)\cap {\cal U}_{\mathfrak{p}}$. Define 
$u_i=\exp^{-1}_{\vp}(\gamma_i)$, $1\le i\le r+1$, and put 
$$\mathcal{L}=m_1u_1+\cdots+m_{r+1}u_{r+1}.$$ 
Set further $M=\max\{|m_1|,\ldots,|m_{r+1}|\}$ and 
$\log W=\max\{\hat h(\gamma_i),1\le i\le r+1\}$.
If $\mathcal{L}\not=0$, then
\begin{align}\label{bound-ultra}
\log |\mathcal{L}|_{\vp}\ge - \kappa_7^{r^2} r^{2r^2} p^{8r^2+28r+23}
& (\log p)^{-3r-4} \, d^{6r+11}\,(\log^+d)^{r+3}\cr
&\times (\log^+ M)\,(\log^+\log W)^{r+3}\,\prod_{i=1}^{r+1}
\max\{1,\hat h(\gamma_i)\}.
\end{align}
\end{corollary}

\medskip
\noindent{\em Proof.} In the following proof, we use the notation of Theorem~\ref{hirata}. Let us begin
by bounding from below the parameter $\mathcal{E}$. Let $n_i$ be the
integer such that $|u_i|_{\vp}=p^{-n_i/e_{\vp}}$. Since $u_i\in {\cal C}_{\vp}$ we have
$|u_i|_{\vp}<p^{-\lambda_{\vp}}$, hence
$n_i/e_{\vp}-\lambda_{\vp}>0$. If $p\not=2$ (hence $\lambda_{\vp}=1/(p-1)$), we have
$$\frac{n_i}{e_{\vp}}-\lambda_{\vp}\ge \frac{1}{(p-1)e_{\vp}}\ge \frac{1}{pe_{\vp}},$$
and if $p=2$ one easily checks that the same bound holds. It follows from this and the definition of $\mathcal{E}$
that we have
\begin{equation}\label{min_E}
\log \mathcal{E} \ge \frac{\log p}{pe_{\vp}}\ge \frac{\log p}{pd}.
\end{equation}
Suppose first that $d\ge\log \mathcal{E}$. Then $\delta=d/\log \mathcal{E}$, and a
rough estimate gives (noticing that $\log(\delta \mathcal{E})\ge 1$ and since $h_E$ is a constant)
\begin{equation}\label{bound_g}
g+\log(\delta \mathcal{E})\ll (\log^+ \log W)\,\log(\delta \mathcal{E}).
\end{equation}
Now, using (\ref{min_E}) we get :
\begin{eqnarray*}
\log(\delta \mathcal{E})&=&\log d + \log \mathcal{E} - \log\log \mathcal{E}\cr
                       &\le& 2\log d +\log p - \log\log p+\log \mathcal{E} \cr
                       &\ll& (\log p)\,\, (\log^+d)\,\, (\log^+ \mathcal{E}).\cr
\end{eqnarray*}
Replacing this estimate in (\ref{bound_g}), we obtain :
\begin{equation*}
g+\log(\delta \mathcal{E})\ll (\log p)\, (\log^+d)\, (\log^+ \log W)\, (\log^+ \mathcal{E}).
\end{equation*}
Using now Hirata's bound (\ref{min_hirata}), we find :
\begin{align*}
\log |\mathcal{L}|_{\vp}\ge - c_1^{r^2} r^{2r^2} & p^{8(r+1)(r+2)}\, (\log p)^{r+3} \, d^{2r+4}\,(\log^+d)^{r+3}\,(\log^+ M)\cr
&\times (\log^+\log W)^{r+3}\,(\log^+\mathcal{E})^{r+3}\,(\log \mathcal{E})^{-4r-7}
\prod_{i=1}^{r+1}\max\{1,\hat h(\gamma_i)\}. 
\end{align*}
Writing finally
$$(\log^+\mathcal{E})^{r+3}\,(\log \mathcal{E})^{-4r-7}=\max\{1,(\log \mathcal{E})^{-1}\}^{r+3}\,(\log \mathcal{E})^{-3r-4}$$
and using the lower bound (\ref{min_E}) to estimate from above this latter
quantity, we obtain (\ref{bound-ultra}) as required.

Suppose now that $d<\log \mathcal{E}$. Then $\delta=1$ and we get in this case
\begin{equation*}
g+\log(\delta \mathcal{E})\ll (\log^+\log W)\, (\log^+ \mathcal{E}).
\end{equation*}
Using (\ref{min_hirata}) and (\ref{min_E}) as before, we find now
\begin{align*}
\log |\mathcal{L}|_{\vp}\ge - c_2^{r^2} r^{2r^2}\, p^{8r^2+26r+19}
&(\log p)^{-2r-3} \, d^{2r+3}\, (\log^+ M)\cr
&\times (\log^+\log W)^{r+3}\,\prod_{i=1}^{r+1}\max\{1,\hat h(\gamma_i)\},
\end{align*}
which again implies the bound (\ref{bound-ultra}). \hfill $\Box$

\medskip 
\noindent{\em Proof of Proposition~\ref{ultra}.} 
Let $Q$, $Q_{r+1}$, $m_1,\ldots,m_r$ and $M$ be as in the Proposition.
We note that since $Q$ is a non-torsion point we have
$M\ge 1$. Denote by $\nu_{\vp}$ the exponent of the group
$E(\Ft_{\mathfrak{p}})/{\cal U}_{\mathfrak{p}}$.
Then $\nu_{\vp}Q,\nu_{\vp}Q_1,\ldots,\nu_{\vp}Q_{r+1}$
belong to $\mathcal{U}_{\vp}$, and the following linear form in
$\vp$-adic elliptic logarithms is well-defined (and non zero
since $\nu_{\vp}Q\not=O$):
$$ \mathcal{L}: = \exp_{\vp}^{-1}(\nu_{\vp}Q) = m_{1}
\exp_{\vp}^{-1}(\nu_{\vp}Q_{1}) + \cdots +
m_r\exp_{\vp}^{-1}(\nu_{\vp}Q_{r})+\exp_{\vp}^{-1}(\nu_{\vp}Q_{r+1}).$$
Since $\varphi:{\cal C}_{\vp}\rightarrow {\cal C}_{\vp}$ is an isometry,
the formula (\ref{exponentielle}) gives
\begin{equation}\label{lambdap} 
|\mathcal{L}|^{-2}_{\vp}=|\varphi(\mathcal{L})|_{\vp}^{-2}= 
|x(\exp_{\vp}(\mathcal{L}))|_{\vp}=|x(\nu_{\vp}Q)|_{\vp}.
\end{equation}
Observe that $|x(Q)|_{\vp}\le |x(\nu_{\vp}Q)|_{\vp}$. Indeed, this is clearly true if $|x(Q)|_{\vp}\le 1$ since $|x(\nu_{\vp}Q)|_{\vp}>1$ by (\ref{U_p}), and this is also true if $|x(Q)|_{\vp}>1$ by Lemma~\ref{filtr}.
This remark together with (\ref{lambdap}) yields
$\log |x(Q)|_{\vp} \leq -2 \log |\mathcal{L}|_{\vp}$.
Applying  Corollary~\ref{min-ultra},
we get an upper bound for $\log |x(Q)|_{\vp}$
involving $\hat{h}(\nu_{\vp}Q_i)$ ($1\le i\le r+1$).
Noticing that $\hat{h}(\nu_{\vp}Q_i)=\nu_{\vp}^2\hat{h}(Q_i)$
and that $\hat{h}(Q_{r+1})=0$, we get
\begin{align*}
\log |x(Q)|_{\vp}\le
2  \kappa_7^{r^2} r^{2r^2} p^{8r^2+28r+23} & \nu_{\vp}^{2r}  (\log p)^{-3r-4} \, d^{6r+11}\,(\log^+d)^{r+3}\cr
&\times (\log^+ M)\,(2\log\nu_{\vp}+\log^+\log V)^{r+3}\,\prod_{i=1}^r\max\{1,\hat h(Q_i)\}.
\end{align*}
But by Lemma~\ref{exposant}, we have
$$ 
\nu_{\vp}\le  p^{\kappa_{5}d}\quad\mbox{hence}\quad
\log\nu_{\vp}\ll d\log p.
$$ 
Proposition~\ref{ultra} follows from these estimates.
\hfill $\Box$

\subsection{Proof of Theorem~\ref{borne-hauteur}}\label{upperbound} 
 
We prove here Theorem~\ref{borne-hauteur}. According to the notation of Section~\ref{notations}, we
write $\mathcal{H}=(<Q_{i}, Q_{j}>)_{1\leq i, j\leq r}$ for the matrix
of the N\'eron-Tate pairing with respect to the chosen basis $(Q_1,\ldots,Q_r)$.

\begin{lemma}\label{neron-tate}
Suppose that $r\ge 1$. Let us denote by $\lambda_{{\rm min}}$ the smallest eigenvalue of the 
matrix ${\cal H}$, and by $\lambda_{{\rm max}}$ its largest eigenvalue.
Let $Q$ be a point of $E(\Ft)$ of the form
$Q=m_1Q_1+\cdots+m_rQ_r+Q_{r+1}$, where $m_1,\ldots,m_r\in\Z$
and $Q_{r+1}$ is a torsion point of $E(\Ft)$. Define further
$M=\max\{|m_1|,\ldots,|m_r|\}$. Then we have
$$\lambda_{{\rm min}}M^2\le \hat{h}(Q) \le r\lambda_{{\rm max}}M^2.$$
\end{lemma}
 
\noindent {\em Proof.} It follows for example from 
\cite[\S~3, inequality~1]{tza} and from its proof. 
\hfill $\Box$ 
 
\begin{lemma}\label{weil-neron-tate} 
For all $Q$ in $E(\Ft)$ we have 
$$\left| \hat{h}(Q) - \frac{1}{2} h_x(Q) \right| \le \kappa_{8}.$$ 
\end{lemma} 
 
\noindent {\em Proof.} 
See \cite[Theorem~1.1]{silverman.difference}. One can take for instance 
$\kappa_{8}=h(\Delta)/12+ h(j(E))/8+1.07$, where 
$\Delta=-16(4A^3+27B^2)$ is the discriminant of the equation 
(\ref{weierstrass}) and $j(E)=-1728(4A)^3/\Delta$ is 
the $j$-invariant of $E$. 
\hfill$\Box$ 
 
\medskip 
\noindent {\em Proof of Theorem~\ref{borne-hauteur}.} 
Let $Q$ be an $S$-integral point of $E(\Ft)$. 
If $\hat h(Q)=0$ then the bound (\ref{bound}) of Theorem~\ref{borne-hauteur} is 
clearly true, since then $h_x(Q)\le 2\kappa_{8}$ by Lemma~\ref{weil-neron-tate}. 
So we will assume in the following that $\hat h(Q)>0$. Thus $Q$ is non-torsion and we have $r\ge 1$.
Write $Q=m_1Q_1+\cdots+m_rQ_r+Q_{r+1}$, where
$m_1,\ldots,m_r$ are integers and $Q_{r+1}$ is a torsion point of $E(\Ft)$.
Define $M:=\max\{|m_1|,\ldots,|m_r|\}$. Applying Proposition~\ref{ultra} to all ultrametric places
$\vp\in S$ and Proposition~\ref{archi} to all archimedean
places $\sigma$, and adding all the inequalities obtained, we get
(recall that $Q$ is $S$-integral, so the places
$v\notin S\cup M_{\Ft}^{\infty}$ do not contribute to the height)
\begin{align*}
h_x(Q) & =\frac{1}{[\Ft:\rat]}
\sum_{v \in S} [\Ft_v:\Q_v] \log
\max \{1, |x(Q)|_{v}\}\cr
&\le c_{1}^{r^2}\, C(E,\Ft)\, (\log^+ M)(\log^+\log M)
\times\frac{1}{[\Ft:\rat]}\bigl(\sum_{v\in S}[\Ft_v:\rat_v]p^{8r^2+\kappa_4 dr}\bigr),\label{qwer}
\end{align*}
where
\begin{equation}\label{const}
C(E,\Ft)=r^{2r^2}d^{9r+14} (\log^+d)^{r+5} (\log^+\!\log V)^{r+5}\, \prod_{i=1}^{r}\max\{1,\hat h(Q_i)\}.
\end{equation}
Now, introducing the set $P(S):=\{p\textrm{ prime}\mid \exists v\in S, v|p\}$, we have :
\begin{align*}
\frac{1}{[\Ft:\rat]}\bigl(\sum_{v\in S}[\Ft_v:\rat_v]p^{8r^2+\kappa_4 dr}\bigr) & \le
\frac{1}{[\Ft:\rat]}\sum_{p\in P(S)}\bigl(\sum_{v|p}[\Ft_v:\rat_v]\bigr) p^{8r^2+\kappa_4 dr}\cr
& = \sum_{p\in P(S)} p^{8r^2+\kappa_4 dr} \leq \prod_{p\in P(S)} p^{8r^2+\kappa_4 dr}  \cr
& = \exp \{(8r^2+\kappa_4 dr) \sum_{p\in P(S)} \log p \}    \leq e^{(8r^2+\kappa_4 dr)\Sigma_S} .
\end{align*}
Hence we deduce
\begin{equation}\label{hhh}
h_x(Q) \le c_{1}^{r^2}\, C(E,\Ft)\, (\log^+M)\, (\log^+\log M)\, e^{(8r^2+\kappa_4 dr)\Sigma_S}.
\end{equation}
Lemma~\ref{weil-neron-tate} yields
$$ \log^+\hat h(Q)\ll \log^+ h_x(Q)\quad
\text{and} \quad \log^+\log\hat h(Q)\ll \log^+\log h_x(Q), $$
and so, by Lemma~\ref{neron-tate} :
\begin{align*}
\log^+ M & \le \frac{1}{2}\bigl(\log^+\hat{h}(Q)+\log^+\lambda_{\rm min}^{-1})
\bigr)\ll (\log^+\! h_x(Q)) . (\log^+\lambda_{\rm min}^{-1})
\end{align*}
and
$$ \log^+\log M\ll (\log^+\log h_x(Q)) .
(\log^+\log \lambda_{\rm min}^{-1}).$$
Substituting these estimates in (\ref{hhh}), we get
\begin{equation}\label{qwertz}
\frac{h_x(Q)}{(\log^+\! h_x(Q))(\log^+\! \log h_x(Q))}\le U
\end{equation}
with
\begin{equation}\label{U}
U = c_{2}^{r^2}\, C(E,\Ft)\, (\log^+\lambda_{\rm min}^{-1})\,
(\log^+\log \lambda_{\rm min}^{-1})\, e^{(8r^2+ \kappa_4 dr)\Sigma_S}.
\end{equation}
We now have
\begin{equation*}
\lambda_{\rm max}\le\trace(\mathcal{H})=\sum_{i=1}^{r}
\hat{h}(Q_{i})\le r \log V,
\end{equation*}
hence
$$ \Reg(E/\Ft)=\det({\mathcal H})\le\lambda_{\rm min}\lambda_{\rm max}^{r-1}\le
\lambda_{\rm min}r^{r-1}(\log V)^{r-1},$$
from which we obtain
$$ \lambda_{\rm min}^{-1}\le\frac{r^{r-1}(\log V)^{r-1}}{\Reg(E/\Ft)}.$$
It follows
$$ \log^+\lambda_{\rm min}^{-1}\ll r(\log^+\! r)\, (\log^+\! \log V)\, (\log^+\!\Reg(E/\Ft)^{-1}) $$
and
$$ \log^+\log \lambda_{\rm min}^{-1}\ll (\log^+\! r)\, (\log^+\!\log^+\!\log V)\, (\log^+\!\log \Reg(E/\Ft)^{-1}). $$
Substituting these inequalities in (\ref{U}) and noticing that
(\ref{qwertz}) implies
$$ h_x(Q)\ll U (\log U) (\log \log U),$$
we get the upper bound (\ref{bound}).
\hfill$\Box$


\bigskip
\noindent \textbf{Remark} :
One would like to bound explicitly the height of the $S$-integral points of $E(K)$ in terms of more manageable objects,
as the set of places $S$, the degree and the discriminant of the number field $K$.
In a forthcoming paper, we show that it is possible to deduce from Theorem~\ref{borne-hauteur} a conditional bound of this kind, 
relying on the conjecture of B. J. Birch and H. P. F. Swinnerton-Dyer \cite{birch-swinnerton-dyer}.
We quote here the result that we obtain.

\begin{prop}\label{conjub}
Let $\Fu$ be a number field, and let $E$ be an elliptic curve given by a Weierstrass equation $y^2=x^3+Ax+B$ with
$A, B\in O_{\Fu}$.  Let $K/\Fu$ be a finite extension, $S$ a finite set of
finite places of $K$, and denote by $d$ the degree $[K:\rat]$ and $D_K$  the absolute value of the discriminant of $K$.
 
Suppose that the $L$-series of $E$ satisfies a  Hasse-Weil functional equation and that the Birch and Swinnerton-Dyer Conjecture holds for $E/K$.

Then, there exist positive numbers $\kappa_{10}$ and  $\kappa_{11}$ (depending on $E/K_0$ only) such that,  for every point $Q$ in $E(O_{K,S})$, we have
$$h_x(Q) \leq \exp\{\kappa_{10}^d +\kappa_{11} \,d^6 (\log^+ D_K)^2\, \left (\Sigma_S + \log (d \log^+ D_K)\right)\}.$$
\end{prop}

Following \cite{bornes}, we deduce from this bound a (weak exponential) inequality of the type of the $abc$-conjecture of D. Masser  and J. Oesterl\'e.


\medskip  
   
\medskip 
   
\bibliography{mabiblio-S-entiers}

 \noindent
\textbf{Vincent Bosser}\\
Laboratoire Nicolas Oresme\\
Universit\'e de Caen\\
F-14032 Caen cedex\\
France

\bigskip

 \noindent
\textbf{Andrea Surroca}\\
Mathematisches Institut\\
Universit\"at Basel\\
Rheinsprung 21\\
CH-4051 Basel\\
Switzerland

\end{document}